\documentclass[reqno]{amsart}

\usepackage{epsfig,amssymb,amsmath,mathrsfs}
\usepackage{hyperref}
\newcommand{\conv}{\operatorname{conv}}

\newcommand{\R}{\mathbb{R}}
\newcommand{\N}{\mathbb{N}}

\newcommand{\union}{\cup}

\newcommand\closure{\operatorname{cl}}

\newcommand\pul{\operatorname{Ls}}
\newcommand\pll{\operatorname{Li}}

\newcommand\indicator{I}
\newcommand\dualspace{V^*}
\newcommand\dualball{{B^\circ}}

\newcommand\dists{\mathcal{D}}
\newcommand\dualdists{\mathcal{{}^* D}}
\newcommand\horocomp{\closure \dists}
\newcommand\dualdbps{{}^*\!\mathcal{A}}
\newcommand\dbps{\mathcal{A}}

\newcommand\dotprod[2]{\langle{#1}|{#2}\rangle}
\newcommand\hnorm{f}

\newtheorem{prop}{Proposition}[section]

\newtheorem{lemma}[prop]{Lemma}

\newtheorem{theorem}[prop]{Theorem}
\newtheorem{example}[prop]{Example}

\begin{document}

\title{The horofunction boundary of finite--dimensional normed spaces}
\date{\today}
\date{\today}

\title{The horofunction boundary of finite-dimensional normed spaces}
\date{\today}
\author{Cormac Walsh}
\address{INRIA, Domaine de Voluceau,
78153 Le Chesnay C\'edex, France}
\email{cormac.walsh@inria.fr}

\renewcommand{\subjclassname}{\textup{2000} Mathematics Subject Classification}

\begin{abstract}
We determine the set of Busemann points of an arbitrary finite--dimensional
normed space. These are the points of the horofunction boundary that
are the limits of ``almost-geodesics''.
We prove that all points in the horofunction boundary are Busemann points
if and only if the set of extreme sets of the dual unit ball is closed
in the Painlev\'e--Kuratowski topology.
\end{abstract}

\maketitle

\section{Introduction}

In~\cite{gromov:hyperbolicmanifolds}, Gromov defines a boundary of a metric
space $(X,d)$ as follows. Let $C(X)$ be the space of continuous real-valued
functions on $X$, with the topology of uniform convergence on compacts
and let $\tilde C(X)$ be the quotient of this space obtained by considering
two functions equivalent if they differ by a constant.
Then, one can use the distance function $d(\cdot,x)$ to inject the space
$X$ into $\tilde C(X)$. If $X$ is proper, meaning that closed balls are
compact, then this map is an embedding.
The topological boundary in $\tilde C(X)$ of the image of this map
is called the \emph{horofunction boundary} of $X$, and its elements are called
\emph{horofunctions}.

Note that this construction is an additive version of the construction of the
Martin boundary in probabilistic potential theory~\cite{dynkin} and of the
Thurston boundary of Teichm\"uller space~\cite{kaimanovich_masur_poissonMGC}.
The analogy with the Martin boundary was developed in~\cite{AGW-m}, where,
in particular, the analogue of the \emph{minimal} Martin boundary was found
to be the set of limits of \emph{almost-geodesics}.
An almost-geodesic, as defined by Rieffel~\cite{rieffel_group},
is a map $\gamma$ from an unbounded set $T\subset \R_+$ containing 0 to $X$,
such that for any $\epsilon>0$,
\begin{equation*}
|d(\gamma(t),\gamma(s))+d(\gamma(s),\gamma(0))-t| < \epsilon
\end{equation*}
for all $t$ and $s$ large enough with $t\ge s$.
Rieffel called the limits of such paths Busemann points.
See~\cite{AGW-m} for a slightly different definition of almost-geodesic
which nevertheless gives rise to the same set of Busemann points.

Rieffel comments that it is an interesting question as to when
all boundary points of a metric space are Busemann points and asks whether
this is the case for general finite--dimensional normed spaces.
We answer this question in the negative and give a necessary and sufficient
criterion for it to be the case.

We exploit the fact that horofunctions in a normed space are convex;
our methods are those of convex geometry, in particular polarity.

Although we have defined horofunctions as equivalence classes of functions,
in the remainder of the paper we only consider their representatives
taking the value $0$ at the origin.

Let $V$ be an arbitrary finite--dimensional normed space with unit ball $B$.
In this paper, all vector spaces are assumed real.
For any extreme set $E$ of the dual unit ball $\dualball$
and point $p$ of $V$, define the function
$\hnorm_{E,p}$ from the dual space $\dualspace$ to $[0,\infty]$ by
\begin{align*}
\hnorm_{E,p}(q):= \indicator_E(q) + \dotprod{q}{p} - \inf_{y\in E}\dotprod{y}{p}
\qquad\text{for all $q\in V^*$}.
\end{align*}
Here $\indicator_E$ is the indicator function, taking value $0$ on $E$
and $+\infty$ everywhere else.

Our first theorem characterises the Busemann points of $V$
as the Legendre--Fenchel transforms of these functions.
\begin{theorem}
\label{thm:theorem1}
The set of Busemann points of a finite--dimensional normed space $(V,||\cdot||)$
is
\begin{equation*}
\{ \hnorm^*_{E,p}  \mid \text{$E$ is a proper extreme set of $\dualball$
and $p\in V$} \}.
\end{equation*}
\end{theorem}

We use this knowledge to characterise those norms for which all horofunctions
are Busemann points.
\begin{theorem}
\label{thm:theorem2}
A necessary and sufficient condition for every horofunction
of a finite--dimensional normed space to be a Busemann point
is that the set of extreme sets of the dual unit ball be closed in the
Painlev\'e--Kuratowski topology.
\end{theorem}

In this paper the symmetry of the norm plays no role; the results hold
equally well for non-symmetric norms.

Karlsson \textit{et.~al.}~determined the horofunction boundary in the case when
the norm is polyhedral~\cite{karl_metz_nosk_horoballs}.
In~\cite{andreev}, Andreev makes a connection between the horofunction boundary
of finite--dimensional normed spaces and flag directed sequences.
The question of when all horofunctions are Busemann points was investigated
in a general setting by Webster and Winchester.
Their paper~\cite{winweb_busemann} contains a criterion for this to be so
when the metric is given by a graph, and~\cite{winweb_metric} contains a
similar criterion for general metric spaces.

\section{Preliminaries}

For a reference on convex analysis, the reader may consult~\cite{beer_book}.

We will use the Painlev\'e--Kuratowski topology on the set of closed
sets of a finite--dimensional normed space $V$.
In this topology, a sequence of closed sets $(C_n)_{n\in\N}$
is said to converge to a closed set $C$
if the upper and lower closed limits of the sequence agree and equal $C$.
These limits are defined to be, respectively,
\begin{align*}
\pul C_n &:= \bigcap_{n\ge 0} \closure \Big( \bigcup_{i>n} C_i \Big)
\qquad
\text{and} \\
\pll C_n &:=
   \bigcap\Big\{ \closure \bigcup_{i\ge 0} C_{n_i}
         \mid \text{$(n_i)_{i\in\N}$ is an increasing sequence in $\N$} \Big\}.
\end{align*}
An alternative characterisation of convergence is that $(C_n)_{n\in\N}$
converges to $C$ if and only if each of the following hold:
\begin{itemize}
\item
for each $x\in C$, there exists $x_n\in C_n$ for $n$ large enough, such that
$(x_n)_n$ converges to $x$.
\item
if $(C_{n_k})_{k\in\N}$ is a subsequence of the sequence of sets and
$x_k\in C_{n_k}$ for each $k\in\N$, then convergence of $(x_k)_{k\in\N}$
to $x$ implies that $x\in C$.
\end{itemize}

The Painlev\'e--Kuratowski topology can be used to define a topology on the
space of lower-semicontinuous functions as follows. Recall that the epigraph
of a function $f$ on $V$ is the set
$\{(x,\alpha)\in V\times\R\mid \alpha\ge f(x)\}$. A sequence of
lower-semicontinuous functions is declared to be convergent in the epigraph
topology if the associated epigraphs converge in the Painlev\'e--Kuratowski
topology on $V\times\R$. For proper metric spaces,
the epigraph topology is identical to the Attouch--Wets topology.

The \emph{Legendre--Fenchel transform} of a function
$f:V \to \R\union\{\infty\}$ is the function $f^*:V^* \to \R\union\{\infty\}$
defined by
\begin{equation*}
f^*(y):= \sup_{x\in V} \big( \dotprod{y}{x} - f(x) \big)
\qquad\text{for all $y\in V^*$}.
\end{equation*}
The Legendre--Fenchel transform is a bijection from the set of
proper lower-semicontinuous convex functions to itself
and is continuous in the epigraph topology.

Let $B:=\{x\in V \mid ||x||\le 1 \}$ be the closed unit ball of the normed
space $(V,||\cdot||)$.
The dual unit ball is the set of linear forms
\begin{equation*}
\dualball
   := \{ y\in V^* \mid \text{$\dotprod{y}{x} \ge -1$ for all $x\in B$} \}.
\end{equation*}
A convex subset $E$ of a convex set $C$ is said to be an \emph{extreme set}
if the endpoints of any line segment in $C$ are contained in $E$ whenever
any interior point of the line segment is.

Let
\begin{align*}
\dualdbps:=\{ \hnorm_{E,p}\mid
             \text{$E$ is an extreme set of $\dualball$ and $p\in V$}\}.
\end{align*}
Note that $\dualdbps$ is precisely the set of functions that are affine on
some extreme set of $\dualball$, take the value $+\infty$ everywhere else,
and have infimum zero.
Let
\begin{align*}
\phi_z(x) := ||z-x||- ||z||
\qquad\text{for all $x$ and $z$ in $V$.}
\end{align*}
Define
\begin{equation*}
\dists := \{ \phi_z \mid z\in V \}
        = \{ ||z-\cdot\,||- ||z|| \mid z\in V \}.
\end{equation*}
Also, let
\begin{align*}
\dbps:=\{ \hnorm^* \mid \hnorm\in\dualdbps \}
\qquad\text{and}\qquad
\dualdists := \{ f^* \mid f\in\dists \}
\end{align*}
be the sets of Legendre--Fenchel transforms of the functions in $\dualdbps$
and $\dists$.

The Legendre--Fenchel transform of $\phi_0(x)=||-x||$ can be
calculated to be $\phi^*_0(y) = \indicator_\dualball(y)$. Also, by expressing
the norm $||\cdot||$ as the transform of its transform, one may arrive at the
formula
\begin{align*}
||z|| = -\inf_{y\in\dualball}\dotprod{y}{z}.
\end{align*}
Using these, one can calculate the transform of
$\phi_z=\phi_0(\,\cdot\,-z)- ||z||$
to be
\begin{align*}
\phi^*_z(y)
   &= \phi^*_0(y) + \dotprod{y}{z} + ||z|| \\
   &= \indicator_\dualball(y) +  \dotprod{y}{z}
             - \inf_{x\in\dualball}\dotprod{x}{z}.
\end{align*}
So $\phi^*_z$ is in $\dualdbps$ for all $z\in V$, and hence $\dists$ is a
subset of $\dbps$.

Denote by $\horocomp$ the horofunction compactification of $V$, that is,
the closure of $\dists$ in the topology of uniform
convergence on compact sets.
As the limit of a sequence of convex 1--Lipschitz functions, each
element of $\horocomp$ is also convex and 1--Lipschitz.
Since the functions in $\dists$ are
equi--Lipschitzian, convergence of a sequence of such functions uniformly
on compacts is equivalent to convergence in the epigraph topology, see
Lemma~7.1.2 and Proposition~7.1.3 of~\cite{beer_book}. So $\horocomp$
is also the closure of $\dists$ in the epigraph topology, which we will
find more convenient to use in the remainder of the paper.

\section{Proof of Theorem~\ref{thm:theorem1}}

The proof of our characterisation of Busemann points will require
a result from~\cite{AGW-m}:
a horofunction is a Busemann point if and only if
it is not the minimum of two 1--Lipschitz functions each different from it.

\begin{lemma}
\label{lem:busemannisinanotd}
Each Busemann point is contained in $\dbps\backslash\dists$.
\end{lemma}
\begin{proof}
Let $g$ be in $(\horocomp)\backslash\dbps$.
So $g$ is convex and 1--Lipschitz.
Its Legendre--Fenchel transform $g^*$ therefore takes the value $+\infty$
outside $B^\circ$.
Since $g^*$ is not in $\dualdbps$, we can therefore find $x$, $y$, and $z$
in $B^\circ$ and $\lambda\in(0,1)$ such that
\begin{align}
\label{eqn:yconv}
y = (1-\lambda)x+\lambda z
\end{align}
and $g^*(y)<(1-\lambda)g^*(x)+\lambda g^*(z)$.
In fact, we can find real numbers $a$ and $b$ such that
$a<g^*(x)$ and $b<g^*(z)$,
and $g^*(y)<(1-\lambda)a+\lambda b$.
It follows from the latter inequality that
\begin{align}
\label{eqn:twostar}
\dotprod{y}{p}-g(p) < (1-\lambda)a+\lambda b
     \qquad\text{for all $p\in V$.}
\end{align}
Define
\begin{align*}
\Pi_1 &:= \Big\{p\in V \mid \dotprod{x}{p}-a \ge g(p) \Big\}
               \qquad\text{and} \\
\Pi_2 &:= \Big\{p\in V \mid \dotprod{z}{p}-b \ge g(p) \Big\}.
\end{align*}
Let $p\in\Pi_1$. Taking the definition of $\Pi_1$ and substituting in the
expression for $x$ obtained from~(\ref{eqn:yconv}) we get
\begin{align*}
\frac{1}{1-\lambda} \dotprod{y}{p}
   - \frac{\lambda}{1-\lambda} \dotprod{z}{p} - a \ge g(p).
\end{align*}
We use the bound on $g(p)$ given by~(\ref{eqn:twostar}) to get, after
some canceling,
\begin{align*}
\dotprod{y}{p} - \dotprod{z}{p} > (1-\lambda)(a-b).
\end{align*}
We use~(\ref{eqn:twostar}) again to deduce that
\begin{align*}
-\dotprod{z}{p} > -g(p)-b.
\end{align*}
Therefore $p$ is not in $\Pi_2$.
We have proved that $\Pi_1$ and $\Pi_2$ have no element in common.

Let
\begin{align*}
g_1 &:= \max(g,\dotprod{x}{\cdot}-a)
               \qquad\text{and} \\
g_2 &:= \max(g,\dotprod{z}{\cdot}-b).
\end{align*}
Each of $g_1$ and $g_2$ are $1$--Lipschitz since
$g$ is $1$--Lipschitz and $x$ and $z$ are in $B^\circ$.
It is immediate from the result of the preceding paragraph that
$g=\min(g_2,g_2)$.

Let $p$ be in the sub-differential of $g^*$ at $x$, which means that
\begin{align*}
\dotprod{q-x}{p} + g^*(x) \le g^*(q)
     \qquad\text{for all $q\in V^*$.}
\end{align*}
Taking the Legendre--Fenchel transform, we get
\begin{align*}
g(s) \le \indicator_{\{p\}}(s)+\dotprod{x}{p} - g^*(x)
     \qquad\text{for all $s\in V$.}
\end{align*}
Now we evaluate at $s=p$ and use the fact that $a<g^*(x)$
to conclude that $g(p) < \dotprod{x}{p} - a$.
Therefore $g_1$ is different from $g$.
One can prove in a similar way that $g_2$ is different from $g$.
Thus we have shown that $g$ is the minimum of two 1-Lipschitz functions
each different from it, and so $g$ cannot be a Busemann point.
\end{proof}

An exposed face of a convex set is the intersection of the set with a
supporting hyperplane.
The following lemma relating this concept to that of extreme set
is probably known, but we can find no reference to it in
the literature.
\begin{lemma}
\label{lem:faceseq}
Let $V$ be a finite--dimensional vector space.
A set $E$ is an extreme set of a convex set $C\subset V$ if and only if
there is a finite sequence of convex sets $F_0,\dots,F_n$ such that
$F_0=C$, $F_n=E$, and $F_{i+1}$ is an exposed face of $F_{i}$ for all
$i\in\{0,\dots,n-1\}$.
\end{lemma}
\begin{proof}
Let $E$ be an extreme set of $C$. If $E$ contains a relative interior point
of $C$, then the extremality of $E$ implies that it equals $C$.
On the contrary, if $E$ is contained entirely within the relative boundary
of $C$, then, since $E$ is convex, by the separation theorem it must
be contained  within an exposed face $F_1$ of $C$.
Since $E\subset F_1 \subset C$ and $E$ is an extreme set of $C$,
it must also be an extreme set of $F_1$.

We may apply the same procedure repeatedly to get the required sequence of sets
$F_0,F_1,\dots,F_n$ such that each is an exposed face of the previous
one and $F_n=E$.

Conversely, assume such a sequence exists. Recall that an exposed face is
an extreme set and that so also is an extreme set of an extreme set.
These two facts imply that $E$ is an extreme set of $C$.
\end{proof}

In the next lemma we will use the following two properties of the epigraph
topology. Firstly, if $(f_n)_{n\in\N}$ is a sequence of proper
lower-semicontinuous
convex functions converging to a limit $f$ such that each takes the value
$+\infty$ outside a fixed bounded region, then $\inf f_n$ converges to
$\inf f$~\cite[Lemma~7.5.3]{beer_book}.
Secondly, if in addition $g$ is a real-valued lower-semicontinuous
convex function that is continuous at a point where $f$ is finite, then
$f_n+g$ converges to $f+g$~\cite[Lemma~7.4.5]{beer_book}.

For each convex subset $C$ of $\dualspace$ and point $p\in V$, define
\begin{align*}
|p|_C := -\inf_{q\in C} \dotprod{q}{p}.
\end{align*}
Observe that $|\cdot|_\dualball = ||\cdot||$.
However, $|\cdot|_C$ will not in general be a norm.
A simple calculation shows that $f^*_{E,p}=|p-\,\cdot\,|_E - |p|_E$
for all extreme sets $E$ of $\dualball$ and points $p$ of $V$.

\begin{lemma}
\label{lem:induct}
Let $F$ be an exposed face of a compact convex set $C$ in $V^*$.
Suppose that there exists a sequence $(p_n)_{n\in\N}$ in $V$
and $\epsilon>0$ such that
\begin{align}
\label{eqn:hyp1}
\sum_{i=0}^{n-1} |p_{i+1}-p_i|_F \le |p_n-p_0|_F + \epsilon
\qquad\text{for all $n\in\N$}
\end{align}
and $|p_{n}-\,\cdot\,|_F - |p_{n}|_F$
converges pointwise to a lower semicontinuous convex function $g$.
Then there exists a sequence $(q_n)_{n\in\N}$ in $V$
and $\epsilon'>0$ such that
\begin{align}
\label{eqn:res1}
\sum_{i=0}^{n-1} |q_{i+1}-q_i|_C \le |q_n-q_0|_C + \epsilon'
\qquad\text{for all $n\in\N$}
\end{align}
and $|q_{n}-\,\cdot\,|_C - |q_{n}|_C$
converges pointwise to $g$.
\end{lemma}
\begin{proof}
There exists an affine function $f$ from $V^*$ to $\R$ that takes the value
zero on $F$ and is positive on $C$. Let $\hat f:= f - f(0)$ be the
linear functional on $V^*$ with the same gradient.

Let $(z_n)_{n\in\N}$ be a sequence of points in $V$ such that
$Z:=\bigcup_n \{z_n\}$ is dense in $V$ and contains the origin.
For each $n\in\N$, define $q_n := p_n + \lambda_n \hat f$,
where the sequence $(\lambda_n)_{n\in\N}$ of reals is chosen so that,
for each $n\in\N$,
\begin{align}
\label{eqn:it1}
|q_{n+1}-q_{n}|_C - |q_{n+1}-q_{n}|_F &< \frac{1}{2^n} \\
\text{and}\qquad\qquad\quad
|q_{n}- z|_C - |q_{n}- z|_F &< \frac{1}{n}
\qquad\text{for all $z\in\{z_0,\dots,z_n\}$}.
\label{eqn:it2}
\end{align}

To see that it is possible to choose the $\lambda_n$ in this way requires
the following argument.
As $\lambda$ tends to infinity, $\lambda f+\indicator_C$ converges in
the epigraph topology to $\indicator_F$.
So $\lambda f+\indicator_C+\dotprod{\cdot}{r}$
converges in the epigraph topology to $\indicator_F+\dotprod{\cdot}{r}$
for any point $r\in V$.
Therefore the infimum of that function converges to the infimum of
$\indicator_F+\dotprod{\cdot}{r}$, which is the same as the infimum
of $\lambda f+\indicator_F+\dotprod{\cdot}{r}$ for any $\lambda\in\R$.
We now use the definitions of $\hat f$, $|\cdot|_C$, and $|\cdot|_F$
to deduce that
\begin{align*}
\lim_{\lambda\to\infty}(|\lambda \hat f + r|_C - |\lambda \hat f + r|_F)
   = 0
\end{align*}
From this and the fact that, for each $n\in\N$,
the set $\{z_0,\dots, z_n\}$ is finite, we see that
(\ref{eqn:it2}) can be satisfied by choosing $\lambda_n$ large enough.
Also, (\ref{eqn:it1}) can be satisfied using large enough $\lambda_{n+1}$
once $\lambda_{n}$ has been fixed.
So one must choose $\lambda_0, \lambda_1, \dots$ in that order.

Because $f$ takes the value $0$ on $F$, we have
\begin{align}
\label{eqn:s2}
|q_{i+1}-q_i|_F &= (\lambda_{i+1}-\lambda_i) f(0) + |p_{i+1}-p_i|_F
\qquad\text{for all $i\in\N$} \\
\text{and}\qquad\qquad
\label{eqn:s3}
|q_{n}-q_0|_F &= (\lambda_{n}-\lambda_0) f(0) + |p_{n}-p_0|_F
\qquad\text{for all $n\in\N$.}
\end{align}
Also, since $F\subset C$,
\begin{align}
\label{eqn:s4}
|q_{n}-q_0|_C \ge |q_{n}-q_0|_F
\qquad\text{for all $n\in\N$.}
\end{align}
Combining~(\ref{eqn:it1}), (\ref{eqn:s4}), (\ref{eqn:s2}), (\ref{eqn:s3}),
and~(\ref{eqn:hyp1}), we get
\begin{align}
\sum_{i=0}^{n-1} |q_{i+1}-q_i|_C - |q_{n}-q_0|_C
   < \sum_{i=0}^{n-1} \frac{1}{2^n} + \epsilon
   \le 2 + \epsilon.
\end{align}
So $(q_n)_{n\in\N}$ satisfies~(\ref{eqn:res1}) with $\epsilon':=2 + \epsilon$.

Let $u\in Z$. For $n\in\N$ large enough,
each of $0$ and $u$ are in $\{z_0,\dots,z_n\}$.
So, from~(\ref{eqn:it2}) and the fact that $|\cdot|_C \ge |\cdot|_F$, we get
\begin{align*}
-\frac{1}{n} \le |q_n - u|_C - |q_n|_C - |q_n - u|_F + |q_n|_F \le \frac{1}{n}
\qquad\text{for all $n$ large enough.}
\end{align*}
But, since $f$ takes the value $0$ on $F$,
we have $|q_n - u|_F - |q_n|_F = |p_n - u|_F - |p_n|_F$,
and by hypothesis this converges to $g(u)$ as $n$ tends to infinity.
We conclude that
$|q_n - u|_C - |q_n|_C$ also converges to $g(u)$ as $n$ tends to infinity.
Note that $|q_n - \,\cdot\,|_C - |q_n|_C$ is $1$--Lipschitz with respect to
any norm on $V$ since its Legendre--Fenchel transform
$\dotprod{\cdot}{q_n}+|q_n|_C+\indicator_C$
takes the value $+\infty$ outside a compact set.
Therefore, the pointwise convergence of this function on
a dense subset of $V$ implies convergence everywhere.
\end{proof}

\begin{lemma}
\label{lem:anotdinbusemann}
Every function in $\dbps\backslash\dists$ is a Busemann point.
\end{lemma}
\begin{proof}
Let $g\in \dbps\backslash\dists$. So $g^*=f_{E,p}$ for some proper extreme
subset  $E$ of $B^\circ$ and point $p$ of $V$.
By Lemma~\ref{lem:faceseq}, there exists a finite sequence of sets
$\dualball=F_0\supset F_1\supset\cdots\supset F_n=E$
such that $F_{i+1}$ is an exposed face of $F_i$ for each
$i\in\{0,\ldots,n-1\}$.

Define the sequence $p_n:=p$ for all $n\in\N$. All the conditions of
Lemma~\ref{lem:induct} are satisfied, taking $F=E$, $C=F_{n-1}$, $\epsilon=0$,
and $g=f^*_{E,p}$. By applying Lemma~\ref{lem:induct} repeatedly,
we arrive at a sequence $(q_n)_{n\in\N}$ in $V$ and $\epsilon'>0$
satisfying the conclusion of this lemma with $C=\dualball$.
So, the sequence of points $(q_n)_{n\in\N}$ converges to the
function $f^*_{E,p}$, which is therefore a horofunction.
Inequality~(\ref{eqn:res1}) says that $(q_n)_{n\in\N}$ is an almost--geodesic
in the sense of~\cite{AGW-m}. Although this notion of almost--geodesic
is slightly different from that of Rieffel, it gives rise to the same set
of Busemann points~\cite[Corollary~7.13]{AGW-m}.
\end{proof}

\begin{proof}[Proof of Theorem~\ref{thm:theorem1}]
This follows from Lemmas~\ref{lem:busemannisinanotd}
and~\ref{lem:anotdinbusemann}.
\end{proof}

\section{Proof of Theorem~\ref{thm:theorem2}}

\begin{lemma}
\label{lem:ae}
If $\dbps$ is closed, then the set of extreme subsets of $\dualball$ is
closed in the Painlev\'e--Kuratowski topology.
\end{lemma}
\begin{proof}
Let $(E_n)_{n\in\N}$ be a sequence of extreme sets of $\dualball$
converging to some set $E$.
The sequence of indicator functions $(\indicator_{E_n})_{n\in\N}$
converges to $\indicator_{E}$. But each
of the functions $\indicator_{E_n};n\in\N$ is in $\dualdbps$,
and so $\indicator_{E}\in\dualdbps$.
It follows that $E$ is an extreme set.
\end{proof}

\begin{lemma}
\label{lem:ea}
If the set of extreme subsets of $\dualball$ is
closed in the Painlev\'e--Kuratowski topology, then $\dbps$ is closed.
\end{lemma}
\begin{proof}
By Lemma~\ref{lem:anotdinbusemann},
we have $\dists\subset\dbps\subset\horocomp$.
So to prove that $\dbps$ is closed, it suffices to prove that the limit $f$ of
any convergent sequence of functions $(f_n)_{n\in\N}$ in $\dualdists$
is in $\dualdbps$.

Let $x$ and $z$ be distinct points in $B^\circ$
and let $y:=(1-\lambda)x + \lambda z$ for some $\lambda\in(0,1)$.
We must show that $f(y)=(1-\lambda)f(x) + \lambda f(z)$.
Since $f$ is convex, there is nothing to prove if $f(y)=\infty$.
We shall therefore assume that $f(y)$ is finite.

Let $|\cdot|$ be any norm on $V^*$.
We claim that there exists a sequence of points $(y_n)_{n\in\N}$ in $B^\circ$
and a constant $\delta>0$ such that the following hold
\begin{itemize}
\item
$(y_n)_{n\in\N}$ converges to $y$,
\item
$f_n(y_n)$ converges to $f(y)$,
\item
for each $n\in \N$, $y_n$ is in some extreme set $E_n$ and is a distance at
least $\delta$ from $\partial E_n$, the relative boundary of $E_n$.
\end{itemize}

Indeed, we know that there exists some sequence $(a_n)_{n\in\N}$ in $B^\circ$
converging to $y$ and satisfying
\begin{equation}
\label{eqn:diamstar}
\lim_{n\to\infty} f_n(a_n)=f(y).
\end{equation}
For each $n\in \N$, let $E_n$ be the smallest extreme set containing $a_n$
and let
\begin{equation}
\delta_n := \inf_{w\in\partial E_n} |a_n-w| .
\end{equation}
If $\delta:= \limsup_{n\to\infty} \delta_n$ is positive,
then our claim holds for some subsequence of $(a_n)_{n\in\N}$,
so assume that $\delta=0$.
For each $n\in \N$, let $b_n\in\partial E_n$ be such that
\begin{equation}
\label{eqn:diamtwo}
|a_n-b_n|=\delta_n.
\end{equation}
Observe that $(b_n)_{n\in\N}$ converges to $y$ as $n\to\infty$.
For each $n\in \N$, let $c_n$ be the point, different from $b_n$, where the line
$a_nb_n$ meets $\partial E_n$.

Define the sequence of points $(s_n)_{n\in\N}$ by
\begin{equation}
\label{eqn:diamone}
s_n := \begin{cases}
   b_n, & \text{if $f_n(b_n) \le f_n(a_n) + \sqrt{\delta_n}$,} \\
   c_n, & \text{otherwise}.
\end{cases}
\end{equation}
Let $(n_i)_{i\in\N}$ be the sequence of $n$ for which the first case
in this definition occurs and let $(m_i)_{i\in\N}$ be the sequence for
which the second occurs.

From equations~(\ref{eqn:diamstar}) and~(\ref{eqn:diamone})
and the fact that $\delta_{n_i}$ tends
to zero as $i\to\infty$, we see that
\begin{align*}
\limsup_{i\to\infty}f_{n_i}(b_{n_i}) \le f(y).
\end{align*}
That the corresponding limit infimum is bounded below by the same quantity
follows immediately from the convergence of $b_n$ to $y$ and the convergence
of $f_n$ to $f$ in the epigraph topology.

Now let $m$ be such that the second case in~(\ref{eqn:diamone}) occurs
for $n=m$. Since $f_m$ is affine on $\dualball$ and $f_m(c_m)\ge 0$, we have
\begin{equation*}
\frac{f_m(a_m)}{|c_m-a_m|}
  \ge \frac{f_m(b_m)-f_m(a_m)}{|b_m-a_m|}
   > \frac{\sqrt{\delta_m}}{\delta_m}
   = \frac{1}{\sqrt{\delta_m}}.
\end{equation*}
Using~(\ref{eqn:diamstar}) and the convergence of $\delta_n$ to zero,
we conclude that
$|c_{m_i}-a_{m_i}|$ tends to zero as $i\to\infty$. So $c_{m_i}$ converges
to $y$. It follows that
\begin{align*}
\liminf_{i\to \infty}f_{m_i}(c_{m_i})\ge f(y).
\end{align*}
We also have that the corresponding limit supremum is no greater than
$f(y)$ because,
for all $i\in\N$, $f_{m_i}(b_{m_i})> f_{m_i}(a_{m_i})$ and the affineness
of $f_{m_i}$ then implies that $f_{m_i}(c_{m_i})\le f_{m_i}(a_{m_i})$.

The results of the previous two paragraphs imply that $(s_n)_{n\in\N}$
converges to $y$ and that $f_n(s_n)$ converges to $f(y)$.
Note that, for all $n\in\N$, $s_n$ is in the boundary of $E_n$, and
hence the smallest extreme set containing $s_n$ has dimension strictly
smaller than that of $E_n$, unless $E_n$ is a singleton. So if we iterate
the procedure constructing $(s_n)_{n\in\N}$, we arrive at a sequence
$(y_n)_{n\in\N}$ which either consists entirely of extreme points or contains
a subsequence satisfying our claim. But the former is impossible
since $\{y\}$ would then be the limit of the extreme sets $\{y_n\};n\in\N$
and therefore, by our assumption, extreme,
contradicting the fact that $y$ is a convex
combination of two points in $\dualball$ distinct from it.
Thus we have proved our claim.

It follows from our claim and the non-negativity of each function $f_n$
that $f_n$ is Lipschitz on $E_n$ with Lipschitz constant $f_n(y_n)/\delta$.
Note that this constant goes to the limit $f(y)/\delta$ as $n\to\infty$.
So we may find a constant $l$ such that
$f_n$ is $l$--Lipschitz on $E_n$ for all $n\in\N$.

Let $F$ be a limit point of the sequence $(E_n)_{n\in\N}$. By assumption,
$F$ must be an extreme set. Taking a subsequence if necessary, we may
assume that $(E_n)_{n\in\N}$ converges to $F$. Since $y_n\in E_n$
for each $n\in\N$ and $(y_n)_{n\in\N}$ converges to $y$, we have $y\in F$.
It follows that each of $x$ and $z$ are in $F$. So there exist sequences
$(x_n)_{n\in\N}$ and $(z_n)_{n\in\N}$ converging to $x$ and $z$ respectively
such that $x_n$ and $z_n$ are in $E_n$ for all $n\in\N$ and $f_n(x_n)$
and $f_n(z_n)$ converge respectively to $f(x)$ and $f(z)$.
Take $y_n':=(1-\lambda)x_n + \lambda z_n$ for all $n\in\N$.
Then $y_n'\in E_n$ for all $n\in\N$ and $y_n'$ converges to $y$.
Since, for each $n\in\N$, the function $f_n$ is Lipschitz on $E_n$
with Lipschitz constant $l$ and $(y_n)_{n\in\N}$ and
$(y_n')_{n\in\N}$ converge to the same limit, $f_n(y'_n)$ must
converge to the same limit as $f_n(y_n)$, which is $f(y)$.
So
\begin{align*}
f(y) &= \lim_{n\to \infty}f_n(y'_n) \\
     &= \lim_{n\to \infty}
             \Big[(1-\lambda)f_n(x_n) + \lambda f_n(z_n) \Big] \\
     &= (1-\lambda)f(x) + \lambda f(z).
\end{align*}
It follows that the set on which $f$ is finite is an extreme set and that $f$
is affine on this set. So $f\in \dualdbps$.
\end{proof}

\begin{lemma}
\label{lem:aclosed}
The set of functions $\dbps$ is closed in the epigraph topology if and only if
the set of extreme subsets of $\dualball$ is closed in the
Painlev\'e--Kuratowski topology.
\end{lemma}
\begin{proof}
This follows from Lemmas~\ref{lem:ae} and~\ref{lem:ea}.
\end{proof}

\begin{proof}[Proof of Theorem~\ref{thm:theorem2}]
Suppose the set of extreme subsets of $\dualball$ is closed.
Then $\dbps$ is closed by Lemma~\ref{lem:aclosed}, and since
$\dists\subset\dbps\subset\horocomp$ by Lemma~\ref{lem:anotdinbusemann},
we have that $\dbps$ equals $\horocomp$.
We conclude that the set of horofunctions
$(\horocomp)\backslash\dists$ equals $\dbps\backslash\dists$,
which by Theorem~\ref{thm:theorem1} is the set of Busemann points.
The converse can be shown by reversing the chain of argument.
\end{proof}

\section{Examples}

\begin{example}
In dimension two, the set of extreme sets of any convex set is always
closed. Therefore, horofunctions of a two-dimensional normed space
are always Busemann points.
\end{example}

\begin{example}
In dimension three, define the norm
\begin{align*}
||(x,y,z)||:= \max\Big( |x|+|z|, \sqrt{x^2+y^2} \Big).
\end{align*}
The unit ball of this norm is
\begin{align*}
B := \Big\{ (x,y,z)\in\R^3
             \mid \text{$|x|+|z|\le 1$ and $x^2+y^2\le 1$} \Big\}.
\end{align*}
The dual unit ball $\dualball$ is the polar of $B$ and is most easily
calculated by
recalling that the polar of an intersection equals the closed convex hull 
of the polars. Thus, $\dualball$ is the convex hull of the square with corners
$(\pm1,0,\pm1)$ and the circle $\{(x,y,z)\mid \text{$x^2+y^2=1$, $z=0$}\}$.

For all $n\in\N$, let $p_n:=(\cos (1/n),\,\sin (1/n),\,0)$.
Observe that the sequence of extreme sets $(\{p_n\})_{n\in\N}$
converges to the set $\{(1,0,0)\}$ as $n\to\infty$.
However, this set is not extreme.

So from Theorem~\ref{thm:theorem2} we would expect the existence of a
horofunction that is not a Busemann point.
Indeed, the function $f:\R^3\to\R,\,(x,y,z)\mapsto -x$ is such a function.

To see it is a horofunction, observe that, for all $n\in\N$,
the sequence of functions $||mp_n-\cdot\,||-||mp_n||$ converges,
as $m$ tends to infinity, to the function
$\xi_n:\R^3\to\R,\,p\mapsto -p_n\cdot p$.
Hence $\xi_n$ is a horofunction for all $n\in\N$, and so $f$ is also
a horofunction since $\xi_n$ converges to $f$ as $n\to\infty$.

To see that $f$ is not a Busemann point, it suffices to write it as the
minimum of the two functions
\begin{align*}
f_1((x,y,z)) &:= \begin{cases}
                  -x+z, & \text{if $z\ge0$}, \\
                  -x,   & \text{otherwise}
                \end{cases} \\
\text{and}\qquad
f_2((x,y,z)) &:= \begin{cases}
                  -x, & \text{if $z\ge0$}, \\
                  -x-z,   & \text{otherwise},
                \end{cases}
\end{align*}
both of which are 1-Lipschitz with respect to the norm $||\cdot||$.
\end{example}

\begin{example}
In dimension three, the set of extreme sets of a convex set is closed
if and only if the set of extreme points is. So, for an example showing that
closure of the set of extreme points of the dual ball is not sufficient for all
horofunctions to be Busemann points, we must go to dimension four.

We define a norm $||\cdot||$ having as dual ball $\dualball$ the closed
convex hull of the four circles
\begin{align*}
S_1^{\pm} &:= \Big\{ (x,y,\pm 1, 0)\in\R^4
             \mid \text{$x^2 + y^2 = 1$} \Big\} \\
S_2^{\pm} &:= \Big\{ (\pm 1, 0, w, z)\in\R^4
             \mid \text{$w^2 + z^2 = 1$} \Big\}.
\end{align*}
The set of extreme points of $\dualball$ is the union of the four circles,
which is closed.

For each $\theta\in(0,\pi/2)$, consider the function $f_\theta:\R^4\to\R$
defined by
\begin{align*}
f_\theta(x,y,w,z):= x \cos\theta + y \sin\theta + z (1-\cos\theta).
\end{align*}
It is easy to show that $f_\theta$ does not take any value greater than $1$
on $\dualball$. Furthermore, for $\theta\in(0,\pi/2)$, $f_\theta$ takes value
$1$ on the set
\begin{align*}
T_\theta := \conv\Big[
              (\cos\theta,\sin\theta,-1,0),
              (\cos\theta,\sin\theta,1,0),
              (1,0,0,1) \Big].
\end{align*}
Therefore this set is a face and hence an extreme set of $\dualball$.
As $\theta$ tends to zero, $T_\theta$ converges to the set
\begin{align*}
\conv\Big[ (1,0,-1,0), (1,0,1,0), (1,0,0,1) \Big].
\end{align*}
However this set is not an extreme set, as can be seen by observing
that it both contains a relative interior point of $\conv S_2^+$ and
is a proper subset of this set.

So the set of extreme sets of $\dualball$ is not closed.

There is enough information in the preceding paragraphs to see that
\begin{align*}
g:\R^4\to\R,\,
(x,y,w,z) \mapsto \max\Big[ x-w, x+w, x+z \Big]
\end{align*}
is a horofunction and that it can be written as the minimum of
two $1$--Lipschitz functions, for example
\begin{align*}
g_1:\R^4\to\R,\, (x,y,w,z) \mapsto
    \max\Big[ x-w, x+w, x+z, x+\frac{w}{\sqrt{2}}+\frac{z}{\sqrt{2}} \Big] \\
\text{and}\qquad
g_2:\R^4\to\R,\, (x,y,w,z) \mapsto
    \max\Big[ x-w, x+w, x+z, x-\frac{w}{\sqrt{2}}+\frac{z}{\sqrt{2}} \Big].
\end{align*}

\end{example}

\begin{example}
It was shown in~\cite{karl_metz_nosk_horoballs} that all horofunctions of a
finite--dimensional normed space with polyhedral norm are Busemann points.
One can recover this result quite easily from Theorem~\ref{thm:theorem2}
by observing that the dual ball in this case is also polyhedral
and has therefore a finite number of extreme sets.
\end{example}

\bibliographystyle{plain}
\bibliography{normed}

\end{document}